\documentclass[a4paper,12pt]{article}

\usepackage{amscd}
\usepackage{amsfonts}
\usepackage{amsmath}
\usepackage{amssymb}
\usepackage{amsthm}
\usepackage[T1]{fontenc} 
\usepackage{here} 
\usepackage{mathrsfs} 
\usepackage{txfonts} 
\usepackage[all]{xy} 
\usepackage{algorithm} 
\usepackage{algpseudocode} 
\usepackage{diagbox} 

\allowdisplaybreaks

\theoremstyle{plain}
\newtheorem{thm}{Theorem}[section]

\newtheorem{conj}[thm]{Conjecture}
\theoremstyle{definition}

\newtheorem{exm}[thm]{Example}

\newcommand{\vs}[1][0.2]{\vspace{#1in}\noindent\ignorespaces}
\newcommand{\ba}{\begin{array*}}
\newcommand{\ea}{\end{array*}}
\newcommand{\be}{\begin{eqnarray*}}
\newcommand{\ee}{\end{eqnarray*}}
\newcommand{\bi}{\begin{itemize}}
\newcommand{\ei}{\end{itemize}}
\newcommand{\bb}{\vs\begin{itembox}}
\newcommand{\eb}{\end{itembox}}
\newcommand{\bc}{\begin{center}}
\newcommand{\ec}{\end{center}}
\newcommand{\bs}{\vs\begin{screen}}
\newcommand{\es}{\end{screen}}

\def\ens#1{{\mathchoice{\left\{ #1 \right\}}{\{ #1 \}}{\{ #1 \}}{\{ #1 \}}}}
\def\set#1#2{{\mathchoice{\left\{ #1 \middle| #2 \right\}}{\{ #1 \mid #2 \}}{\{ #1 \mid #2 \}}{\{ #1 \mid #2 \}}}}
\def\r#1{\text{\rm #1}}

\newcommand{\bF}{\mathbb{F}}

\newcommand{\bN}{\mathbb{N}}

\newcommand{\bP}{\mathbb{P}}
\newcommand{\bQ}{\mathbb{Q}}

\newcommand{\bZ}{\mathbb{Z}}

\newcommand{\cA}{\mathscr{A}}

\newcommand{\cZ}{\mathscr{Z}}

\newcommand{\rL}{\r{L}}

\newcommand{\rN}{\r{N}}

\newcommand{\rR}{\r{R}}

\newcommand{\F}{\bF}
\newcommand{\N}{\bN}
\newcommand{\Q}{\bQ}

\newcommand{\Z}{\bZ}

\newcommand{\Fp}{\mathbb{F}_p}

\newcommand{\gen}{\r{generated}}

\newcommand{\rdp}{\r{dp}}

\algnewcommand\algorithmicbreak{{\bf break}}
\algnewcommand\Break{\algorithmicbreak{}}
\algnewcommand\algorithmiccontinue{{\bf continue}}
\algnewcommand\Continue{\algorithmiccontinue{}}

\title{Bounded Additive Relation and Application to Finite Multiple Zeta Values}
\author{Tomoki Mihara}
\date{}

\begin{document}

\maketitle
\begin{abstract}
We formulate an algebraic problem to find a generating system of a finite subset of an Abelian group with respect to linear relations whose coefficients are bounded by a constant, and recall MITM algorithm for the problem. As an application of MITM algorithm for the Abelian group
\be
\Z/106700590455862347842907841856033238416352421 \Z
\ee
combined with Chinese remainder algorithm, we give a table of expected linear relations of finite multiple zeta values of weight $10$.
\end{abstract}

\tableofcontents

\section{Introduction}
\label{Introduction}

A finite multiple zeta value (FMZV) is a finite analogue of a multiple zeta value introduced by D.\ Zagier. Although no FMZV is known to be non-zero, there are many conjectures on the non-triviality of specific values. Especially, the following conjecture (cf.\ \cite{Kan17} p. 177 and \cite{KZ} p.\ 2) is one of main topics in the study of FMZVs:

\begin{conj}[Dimension conjecture]
For a $w \in \N$, let $d_w$ denote the dimension of the $\Q$-vector space generated by FMZVs of weight $w$. Then the sequence $(d_w)_{w=0}^{\infty}$ satisfies the following recursive relation:
\be
d_w =
\left\{
\begin{array}{ll}
0 & (w \in \ens{1,2}) \\
1 & (w \in \ens{0,3}) \\
d_{w-2} + d_{w-3} & (w > 3)
\end{array}
\right.
\ee
\end{conj}

According to \cite{Kan17} p.\ 177, the dimension conjecture was originally obtained from numerical analysis in \cite{KZ} based on a non-trivial method by D.\ Zagier. However, \cite{KZ} was unpublished for years, and other experts tend not to write down explicit algorithms for numerical analysis of FMZVs.

\vs
Although we personally invented a method for the numerical analysis of FMZVs, we were not certain whether it had novelty or not because none of algorithms by experts was available to us, and hence we did not consider to make it public. However, recently \cite{KZ} has been publicly available as a preprint, and we noticed that our method seems to have novelty to some extent in the sense that the unpublished method is included in our method. Therefore, we decided to start writing this paper to explain our method for the numerical analysis.

\vs
We briefly explain contents of this paper. In \ref{Convention}, we introduce convention for this paper. In \ref{Bounded Additive Relation Problem}, we introduce algebraic problems called a {\it bounded additive relation problem} and a {\it minimal generating system}. In \ref{MITM Algorithm}, we recall a classical method called {\it MITM algorithm}, and explain how to solve a bounded additive relation problem. In \ref{Dynamic MITM Algorithm}, we introduce a dynamic variant of MITM algorithm, and explain how to solve a minimal generating system. In \ref{Numerical Analysis of Finite Multiple Zeta Values}, we briefly recall the definition of FMZVs, recall dynamic programming to compute a mod $p$ harmonic sum in \cite{KZ}, introduce a new parallel computation method based on depth first search improving the dynamic programming, explain a new method for numerical analysis of FMZVs based on combination of MITM algorithm and another classical method called {\it Chinese remainder algorithm}, and exhibits a table of expected linear relations of FMZVs of weight $10$.

\section{Convention}
\label{Convention}

We denote by $\N$ the set of non-negative integers. For a $d \in \Z$ and a binary relation $R$ on $\Z$, we set $\N_{R d} \coloneqq \set{i \in \N}{i R d}$. For sets $X$ and $Y$, we denote by $X^Y$ the set of maps $Y \to X$. We note that every $d \in \N$ is identified with $\N_{< d}$ in set theory, and hence $X^d$ formally means $X^{\N_{< d}}$, which is naturally identified with the set of $d$-tuples in $X$.

\vs
When we write a pseudocode, a for-loop along a subset of $\N$ denotes the loop of the ascending order, and a for-loop along a general set $S$ denotes a loop in an arbitrary order.

\section{Bounded Additive Relation Problem}
\label{Bounded Additive Relation Problem}

Let $M$ be an Abelian group, $(D,B) \in \N^2$, $\vec{x} = (x_d)_{d=0}^{D-1} \in M^D$, $\vec{c} = (c_b)_{b=0}^{B-1} \in \Z^B$, and $S \subset \N_{< B}^D$. Consider the problem to detect the existence of a $\vec{b} = (b_d)_{d=0}^{D-1} \in S$ with $\sum_{d=0}^{D-1} c_{b_d} x_d = 0$ and find such a $\vec{b}$ if exists, which we call {\it the bounded additive relation problem for $(M,D,\vec{x},B,\vec{c},S)$}. We will explain MITM algorithm to solve a bounded additive relation problem in \S \ref{MITM Algorithm}.

\begin{exm}
Let $(N,D,B) \in \N^3$, and $\vec{x} = (x_d)_{d=0}^{D-1} \in (\Z/N \Z)^D$. Set $C \coloneqq (-B+b)_{b=0}^{2B}$, and $S \coloneqq \N_{< B}^D \setminus (b)_{d=0}^{D-1}$. Then the bounded additive relation problem for $(\Z/N \Z,D,\vec{x},2B+1,C,S)$ is equivalent to the problem to detect the existence of a $\vec{c} = (c_d)_{d=0}^{D-1} \in \Z^D$ satisfying the following and find such a $\vec{c}$ if exists:
\bi
\item[(1)] There exists a $d \in \N_{< D}$ with $c_d \neq 0$.
\item[(2)] For any $d \in \N_{< D}$, the inequality $-B \leq c_d \leq B$ holds.
\item[(3)] The equality $\sum_{d=0}^{D-1} c_d x_d = 0$ holds.
\ei
For example, if $(N,D,\vec{x},B) = (7,2,(2 + 7 \Z, 3 + 7 \Z),2)$, then $\vec{c} = (2,1)$ is a solution. If $(N,D,\vec{x},B) = (100,2,(2 + 7 \Z, 3 + 7 \Z),3)$, then $\vec{c} = (-3,2)$ is a solution. On the other hand, if $(N,D,\vec{x},B) = (100,2,(2 + 7 \Z, 3 + 7 \Z),2)$, then there is no solution.
\end{exm}

For a $y \in M$, we say that $\vec{x}$ {\it generates $y$ over $\vec{c}$} if there exists a $(b_d)_{d=0}^{D} \in \N_{< B}^{D+1}$ satisfying the following:
\bi
\item[(1)] The inequality $c_{b_D} \neq 0$ holds.
\item[(2)] The equality $(\sum_{d=0}^{D-1} c_{b_d} x_d) + c_{b_D} y = 0$ holds.
\ei
The problem to detect whether $\vec{x}$ generates a given $x \in M$ over $\vec{c}$ or not is equivalent to the bounded additive relation problem for $(M,D+1,\vec{x} \frown x,B,\vec{c},\set{(b_d)_{d=0}^{D} \in \N_{< B}^{D+1}}{c_{b_D} \neq 0})$, where $\frown$ denotes the concatenation.

\vs
We say that $\vec{x}$ {\it generates $M$ over $\vec{c}$} if it generates any $x \in M$ over $\vec{c}$. When $M$ is a $\mathbb{Q}$-vector space, then the condition is quite similar to the definition of a set of generators. However, the notion makes sense even when $M$ is a torsion Abelian group. Since the problem to detect whether $\vec{x}$ generates a given $x \in M$ over $\vec{c}$ or not is reduced to a bounded additive relation problem, the problem to detect whether $\vec{x}$ generates $M$ over $\vec{c}$ or not is reduced to bounded additive relation problems.

\vs
Let $S \subset M$ be a finite subset. We say that $\vec{x}$ is a {\it minimal generating system of $S$ over $\vec{c}$} if it generates $S$ over $\vec{c}$ and there is no $D' \in \N_{< D}$ such that $(x_d)_{d=0}^{D'-1}$ generates $x_{D'}$ over $\vec{x}$. Needless to say, this notion is a bounded analogue of that of a linear basis. Consider the problem to construct a minimal generating system of $S$ over $\vec{c}$, which we call {\it the minimal generating system problem for $(S,\vec{c})$}.

\vs
In order to solve the minimal generating system problem for $(S,\vec{c})$, it suffices to apply a greedy method. First, initialise $\vec{x}$ as the empty array. Search all $x \in M$ in any order, and if $\vec{x}$ does not generate an $x \in S$ over $\vec{c}$, then append $x$ to $\vec{x}$. Then the resulting $\vec{x}$ is a minimal generating system of $S$ over $\vec{c}$ by definition. In particular, the construction problem is reduced to bounded additive relation problems. 

\section{MITM Algorithm}
\label{MITM Algorithm}

Let $X$ be a finite set, and $Y$ be a set, $f$ a map $X \to Y$, and $y \in Y$. Consider the problem to detect the existence of an $x \in X$ with $f(x) = y$ and find such an $x$ if exists, which we call {\it the decipher problem for $(X,Y,y,f)$}.

\vs
A brute-force on $X$ requires $\# X$ time computations of $f$ in a worst case. The meet-in-the-middle (MITM) Algorithm is a space-time trade off method to effectively solve the problem in the case where $f$ is a composite map. Although MITM algorithm is applicable to various problems, we recall its formulation specialised to a problem based on a binary operation.

\vs
Let $X_0$, $X_1$, and $Y'$ be finite sets such that $X$ is a subset of $X_0 \times X_1$, $g$ a map $Y' \times X_1 \to Y$, $g'$ a map $X_1 \to Y'$ such that $g(y',x_1) = y$ if and only if $y' = g'(x_1)$ for any $(y',x_1) \in Y' \times X_1$, and $h$ a map $X_0 \to Y'$ such that $f(x_0,x_1) = g(h(x_0),x_1)$ for any $(x_0,x_1) \in X$. Then we have
\be
f^{-1}(y) & = & \bigcup_{x_1 \in X_1} \set{x_0 \in X_0}{(x_0,x_1) \in X \land f(x_0,x_1) = y} \times \ens{x_1} \\
& = & \bigcup_{x_1 \in X_1} \set{x_0 \in X_0}{(x_0,x_1) \in X \land g(h(x_0),x_1) = y} \times \ens{x_1} \\
& = & \bigcup_{x_1 \in X_1} \set{x_0 \in X_0}{(x_0,x_1) \in X \land h(x_0) = g'(x_1)} \times \ens{x_1} \\
& = & \bigcup_{x_1 \in X_1} \set{x_0 \in h^{-1}(g'(x_1))}{(x_0,x_1) \in X} \times \ens{x_1}.
\ee
MITM algorithm first computes $h^{-1}(y')$ for each $y' \in Y'$, and next search an $x_1 \in X_1$ such that there is an $x_0 \in h^{-1}(g'(x_1))$ such that $(x_0,x_1) \in X$. Here is a pseudocode for MITM algorithm:

\begin{figure}[H]
\begin{algorithm}[H]
\caption{MITM algorithm}
\label{MITM algorithm}
\begin{algorithmic}[1]
\Function {MITM}{$X,X_0,X_1,Y',Y,g,g',h,y$}
	\State $\vec{\Sigma} \coloneqq (\Sigma_{y'})_{y' \in Y'} \gets$ The associative array of empty arrays indexed by $Y'$
	\ForAll {$x_0 \in X_0$}
		\State Append $\Sigma_{h(x_0)}$ to $x_0$
	\EndFor
	\ForAll {$x_1 \in X_1$}
		\ForAll {$x_0 \in \Sigma_{g'(x_1)}$}
			\If {$(x_0,x_1) \in X$}
				\State \Return $(x_0,x_1)$
			\EndIf
		\EndFor
	\EndFor
	\State \Return False \Comment{flag to indicate the non-existence of a solution}
\EndFunction
\end{algorithmic}
\end{algorithm}
\end{figure}

We call $\vec{\Sigma}$ {\it the MITM dictionary}. The construction of the MITM dictionary requires $\# X_0$ time computations of $h$ and appending an element to an array. Although the computational complexity of the construction of $(\Sigma_{y'})_{y' \in Y}$ depends on implementation of an associative array, the time complexity the space complexity typically depends on $\# X_0$ in the linear order. This is the reason why MITM algorithm is a space-time trade off method.

\vs
The searching process requires $\# X_1$ times computation of $g'$ and $X_0$ times computation of the membership relation for $X$.

\vs
In particular, if $\# X_0 + \# X_1$ is significantly smaller than $\# X$, then MITM algorithm performs better than a brute-force on $X$.

\vs
Let $M$ be an Abelian group, $(D,B) \in \N^2$, $\vec{x} = (x_d)_{d=0}^{D-1} \in M^D$, $\vec{c} = (c_b)_{b=0}^{B-1} \in \Z^B$, and $S \subset \N_{< B}^D$. The bounded additive relation problem for $(M,D,\vec{x},B,\vec{c},S)$ is the same as the decipher problem for $(S,M,0,f)$, where $0$ is the zero element of $M$ and $f$ is the map
\be
S & \to & M \\
(b_d)_{d=0}^{D-1} & \mapsto & \sum_{d=0}^{D-1} c_{b_d} x_d.
\ee
Suppose $D > 1$ and $\# M < \infty$. Set $D^{\rL} \coloneqq \lfloor D/2 \rfloor$, $D^{\rR} \coloneqq D - D^{\rL}$ $X_0 \coloneqq \N_{< B}^{D^{\rL}}$, $X_1 \coloneqq \N_{< B}^{D^{\rR}}$, and $Y' \coloneqq M$. We define
\be
g \colon Y' \times X_1 & \to & M \\
(y',(i_d)_{d=0}^{D^{\rR}-1}) & \mapsto & y' + \sum_{d=D^{\rL}}^{D-1} c_{b_{d-D^{\rL}}} x_d \\
g' \colon X_1 & \to & Y' \\
(i_d)_{d=0}^{D^{\rR}-1} & \mapsto & - \sum_{d=D^{\rL}}^{D-1} c_{b_{d-D^{\rL}}} x_d \\
h \colon X_0 & \to & Y' \\
(i_d)_{d=0}^{D^{\rL}-1} & \mapsto & \sum_{d=0}^{D^{\rL}-1} c_{b_d} x_d
\ee
Then $X$ is identified with a subset of $X_0 \times X_1$ through a natural bijection $X_0 \times X_1 \cong \N_{< B}^D$, and $(X_0,X_1,Y',g,g',h)$ satisfies the condition for MITM algorithm for the decipher problem for $(S,M,0,f)$.

\section{Dynamic MITM Algorithm}
\label{Dynamic MITM Algorithm}

We have explained a greedy construction of a minimal generating system by solving bounded additive relation problems, and MITM algorithm to solve a  bounded additive relation problem. Naively combining them, we obtain an algorithm to construct a minimal generating system.

\vs
However, the naive method is redundant in the sense that the construction of an MITM dictionary is repeated. Instead, it is better to store it unless we need to actually reconstruct it. For example, we do not need to reconstruct an MITM dictionary if $\vec{x}$ is not updated. When $\vec{x}$ is updated, $D_{\rL}$ can be updated to $D^{\rL} + 1$. If updating $D^{\rL}$ reduces expected time complexity of the rest process, then it is good to update $D^{\rL}$. Here are pseudocodes for the improved process:

\begin{figure}[H]
\begin{algorithm}[H]
\caption{Expected time complexity for $D^{\rR}$}
\label{Expected time complexity}
\begin{algorithmic}[1]
\Function {CostD}{$B,D^{\rR},H,h$}
	\State \Return $B^{D^{\rR}}(H-h)$
\EndFunction
\end{algorithmic}
\end{algorithm}
\end{figure}

\begin{figure}[H]
\begin{algorithm}[H]
\caption{Dynamic MITM algorithm to compute a minimal generator system of $S \subset M$ over $\vec{c}$}
\label{Dynamic MITM algorithm}
\begin{algorithmic}[1]
\Function {DynamicMITM}{$M,B,\vec{c}$}
	\State $D \gets 0$
	\State $(D^{\rL},D^{\rR}) \gets (0,0)$
	\State $\vec{x} \gets ()$ \Comment{$x_d$ for a $d \in \N$ smaller than the current length of $\vec{x}$ denotes the $(1+d)$-th entry of $\vec{x}$}
	\State $\vec{\Sigma} \coloneqq (\Sigma_{y'})_{y' \in M} \gets$ the associative array of empty arrays indexed by $M$
	\State $H \gets \# S$
	\State $\vec{s} = (s_h)_{h=0}^{H-1} \gets$ an enumeration of $S$
	\ForAll {$h \in \N_{< H}$}
		\State $\gen \gets$ False
		\ForAll {$b \in \N_{< B}$}
			\If {$c_b = 0$}
				\State \Continue
			\EndIf
			\ForAll {$(b^{\rR}_d)_{d=0}^{D^{\rR}-1} \in \N_{< B}^{D_{\rN}}$}
				\State $y \gets (\sum_{d=0}^{D^{\rR}-1} c_{b^{\rR}_d} x_{D_{\rL}+d}) + c_b s_h$
				\If {$\Sigma_{-y}$ is empty}
					\State \Continue
				\EndIf
				\State $\gen \gets$ True
				\State \Break
			\EndFor
			\If {$\gen$}
				\State \Break
			\EndIf
		\EndFor
		\If {$\neg \gen$}:
			\State Append $s_h$ to $\vec{x}$
			\State $D \gets D + 1$
			\State $w_0 \gets$ \Call{CostD}{$B,D^{\rR},H,h$}
			\State $w_1 \gets$ \Call{CostD}{$B,D^{\rR}+1,H,h$}
\algstore{Dynamic MITM algorithm 1}
\end{algorithmic}
\end{algorithm}
\end{figure}
\addtocounter{algorithm}{-1}
\begin{figure}[H]
\begin{algorithm}[H]
\begin{algorithmic}[1]
\algrestore{Dynamic MITM algorithm 1}
			\If {$w_0 > B^{D^{\rL}+1} + w_1$} \Comment{optimisation of expected computational complexity}
				\State $D^{\rL} \gets D^{\rL} + 1$
				\State $\vec{\Sigma} \gets$ the associative array of empty arrays indexed by $M$
				\ForAll {$\vec{b}^{\rL} \coloneqq (b_d)_{d=0}^{D^{\rL}-1} \in \N_{< B}^{D^{\rL}}$}
					\State $y \gets \sum_{d=0}^{D^{\rL}-1} c_{b_d} x_d$
					\State Append $\vec{b}^{\rL}$ to $\Sigma_y$
				\EndFor
			\Else
				\State $D^{\rR} \gets D^{\rR} + 1$
			\EndIf
		\EndIf
	\EndFor
	\State \Return $\vec{x}$
\EndFunction
\end{algorithmic}
\end{algorithm}
\end{figure}

We note that if we are not interested in actual bounded additive relations, it suffices to replace $\vec{\Sigma}$ by a set of keys rather than an associative array, because we only use whether $\Sigma_y$ is empty or not for a $y \in M$.

\section{Numerical Analysis of Finite Multiple Zeta Values}
\label{Numerical Analysis of Finite Multiple Zeta Values}

Let $\bP$ denote the set of prime numbers, and set
\be
\cA \coloneqq \prod_{p \in \bP} \Fp \bigg/ \bigoplus_{p \in \bP} \Fp.
\ee
for an $L \in \N$ and a $\vec{k} = (k_{\ell})_{\ell=0}^{L-1} \in \N_{> 0}^L$, the finite multiple zeta value (FMZV) of index $\vec{k}$ is the residue class $\zeta^{\cA}(\vec{k})$ of the sequence
\be
\left( \sum_{\vec{m} = (m_{\ell})_{\ell=0}^{L-1}} \frac{1}{\prod_{\ell=0}^{L-1} m_{\ell}^{k_{\ell}}} \right)_{p \in \bP}
\ee
of ``mod $p$ multiple harmonic sums'' in $\cA$, where $\vec{m}$ runs through all strictly increasing sequences of positive integers smaller than $p$. We note that $\zeta^{\cA}(()) = 1$ by definition.

\vs
Since $\cA$ forms a $\Q$-algebra, it is natural to ask the $\Q$-linear dimension of a $\Q$-vector space generated by finitely many FMZVs. In particular, for a $w \in \N$, we set
\be
K_w \coloneqq \set{(k_{\ell})_{\ell=0}^{L-1} \in \N_{> 0}^L}{L \in \N \land \sum_{\ell=0}^{L-1} k_{\ell} = w} \subset \bigcup_{l \in \N} \N_{> 0}^L,
\ee
denote by $\cZ_{\cA,w} \subset \cA$ the $\Q$-linear subspace generated by $\set{\zeta^{\cA}(\vec{k})}{\vec{k} \in K_w}$, and put $d_w \coloneqq \dim_{\Q} \cZ_{\cA,w}$.

\subsection{Dynamic Programming}
\label{Dynamic Programming}

First, we recall dynamic programming (DP) to compute a mod $p$ multiple harmonic sum explained in \cite{KZ} p.\ 21. For an $L \in \N$, $\vec{k} = (k_{\ell})_{\ell=0}^{L-1} \in \N_{> 0}^L$, $p \in \bP$, and $j \in \N_{< p}$, we set
\be
\rdp_p(\vec{k},j) \coloneqq \sum_{\vec{m} = (m_{\ell})_{\ell=0}^{L-1}} \frac{1}{\prod_{\ell=0}^{L-1} m_{\ell}^{_{\ell}}},
\ee
where $\vec{m}$ runs through all strictly increasing sequences of positive integers smaller than $p$ with $m_{L-1} \leq j$ if $L > 0$. By the definition, $\zeta^{\cA}(\vec{k})$ is the residue class of $(\rdp_p(\vec{k},p-1))_{p \in \bP}$ for any $\vec{k} \in \N_{> 0}^L$ with $L \in \N$.

\vs
Let $p \in \bP$. As the base case, we have
\be
\rdp_p((),j) = 1
\ee
for any $j \in \N_{< p}$. For any $(\vec{k},k,j) \in \N_{> 0}^L \times \N \times \N_{< p}$ with $L \in \N$, we have
\be
\rdp_p(\vec{k} \frown k,j) = 
\left\{
\begin{array}{ll}
0 & (j = 0) \\
\rdp_p(\vec{k} \frown k,j - 1) + \rdp_p(\vec{k},j) \times \frac{1}{j^k} & (j > 0)
\end{array}
\right..
\ee
As a consequence, the table $\rdp_p$ for prefixes of $\vec{k}$ is computed in $2$-dimensional DP. Practically speaking, the $2$-dimensional DP does not require $2$-dimensional memory for all intermediate values, and it suffices to keep the last row (resp.\ column) and the current row (resp.\ column) to reduce space complexity.

\vs
In particular, there are two directions to proceed the $2$-dimensional DP: One (horizontal DP) is the DP restoring intermediate values for all prefixes of $\vec{k}$ and only the last and the current value of $j$, and the other one (vertical DP) is the DP restoring intermediate values for all $j$'s and only the last and the current prefixes of $\vec{k}$.

\vs
In horizontal DP, the update for the current value of $j$ refers to the last value of $j$ with a shorter index. In vertical DP, the update for the current prefix refers to the last prefix with a smaller value of $j$. Therefore, both DP can be executed using an in-place array. Although they share $\Theta(pL)$ time complexity, horizontal DP is of $\Theta(L)$ space complexity (using binary modular method or Euclidean method for computing modular inverses), and vertical DP is of $\Theta(p)$ space complexity (using a recursive method for computing modular inverses). Therefore, it is good to choose the direction in a way depending on $L$ and $p$.

\vs
Next, we explain parallel computation of mod $p$ multiple harmonic sums as improvement of the method above. A {\it tree of indices} is a rooted tree $T$ whose node is a finite sequence of positive integers such that for any node $\vec{k}$ of $T$, if $\vec{k}$ is non-empty, then the array given by removing the last entry of $\vec{k}$ is the parent node of $\vec{k}$.

\begin{exm}
\label{tree of indices}
\bi
\item[(1)] For any finite sequence $\vec{k}$ of positive integers, the set of prefixes of $\vec{k}$ admits a unique structure of a tree of indices. Say, if $\vec{k} = (1,3,2,4)$, then the corresponding tree of indices is
\be
() \to (1) \to (1,3) \to (1,3,2) \to (1,3,2,4).
\ee
The construction of the linear graph by handling the length of a prefix instead of copying a prefix is executed in $\Theta(k)$ time and space complexity.
\item[(2)] For any $w \in \N_{< 0}$, the set $K_{\leq w} \coloneqq \bigcup_{w' \in \N_{\leq w}} K_{w'}$ of finite sequence of positive integers whose sum is not greater than $w$ admits a unique structure of a tree of indices. The construction of $K_{\leq w}$ by depth first search is executed in $\Theta(2^w)$ time and space complexity.
\ei
\end{exm}

Let $T$ be a tree of indices. We denote by $V_T$ the set of vertices of $T$. Since $T$ is an acyclic directed graph compatible with the recursive relation of $\rdp_p$, horizontal DP works also for $T$. Namely, it suffices to proceed the $2$-dimensional DP by repeating depth first search on $T$ for each $j \in \N_{< p}$ and restoring intermediate values for all nodes of $T$ and only the last and the current value of $j$. Here is a pseudocode of the process:

\begin{figure}[H]
\begin{algorithm}[H]
\caption{Parallel computation of the mod $p$ harmonic sums for a tree $T$ of indices}
\label{parallel horizontal DP}
\begin{algorithmic}[1]
\Function {ParallelHorizontalDP}{$p,T$}
	\State $\rdp \coloneqq (\rdp_{\vec{k}})_{\vec{k} \in V_T} \gets (0)_{\vec{k} \in V_T}$
	\State $\rdp_{()} \gets 1$
	\Function {DFS}{$\vec{k},j$}
		\ForAll {child nodes $\vec{k}_{+}$ of $\vec{k}$ in $T$}
			\State Execute \Call{DFS}{$\vec{k}_{+},j$}
			\State $k \gets$ the last entry of $\vec{k}_{+}$
			\State $\rdp_{\vec{k}_{+}} \gets (\rdp_{\vec{k}_{+}} + \rdp_{\vec{k}} \times \frac{1}{j^k}) \bmod p$
		\EndFor
	\EndFunction
	\ForAll {$j \in \N_{< p - 1}$}
		\State Execute \Call{DFS}{$(),j+1$}
	\EndFor
	\State \Return $\rdp$
\EndFunction
\end{algorithmic}
\end{algorithm}
\end{figure}

We note that Algorithm \ref{parallel horizontal DP} applied to a tree of indices in Example \ref{tree of indices} (1) is equivalent to the DP in \cite{KZ}, although we are not certain whether it is the naive $2$-dimensional DP, horizontal/vertical DP with an in-place array, or some other implementation.

\vs
If we naively repeat horizontal DP to indices in $V_T$, then the whole process is executed in $\Theta(p W_T)$ time complexity and $\Theta(\# V_T)$ space complexity, where $W_T$ denotes the sum of the depth of nodes in $T$. On the other hand, the parallel computation is executed in $\Theta(p \# V_T)$ time complexity and $\Theta(\# V_T)$ space complexity, and hence actually improves the naive method.

\subsection{Chinese Remainder Algorithm}
\label{Chinese Remainder Algorithm}

Let $w \in N$. Take a subset $S \subset K_w$ of cardinality $d_w$ such that $(\zeta^{\cA}(\vec{k}))_{\vec{k} \in S}$ forms a $\Q$-linear basis of $\cZ_{\cA,w}$. For a $B \in \N$, a rational number is said to be {\it $B$-expressible} if it admits a fractional expression by integers in $\Z \cap [-B,B]$, and we denote by $Q_B \subset \Q$ the subset of $B$-expressible rational numbers.

\vs
Since $K_w$ is finite, there exists a $B \in \N$ such that for any $\vec{k} \in K_w$, the coefficients of the presentation of $\zeta^{\cA}(\vec{k})$ by $(\zeta^{\cA}(\vec{k}))_{\vec{k} \in S}$ are $B$-expressible. By the finiteness of $K_w$ again, we may replace $B$ by a bigger constant so that the same presentation holds for any mod $p$ harmonic sum for any $p \in \bP_{> B}$.

\vs
On the other hand, even if we conversely have a similar expression of a mod $p$ harmonic sum by $B$-expressible rational numbers for a $p \in \bP_{> B}$, it does not mean that the same expression holds for the corresponding FMZV. Namely, a minimal generating system for the pair of the $\Fp$-vector space of mod $p$ harmonic sums of weight $w$ and an array $\vec{c}$ enumerating $\Z \cap [-B,B]$ does not necessarily lift to a $\Q$-linear basis of $\cZ_{\cA,w}$.

\vs
Especially when $p$ is not significantly greater than $(\# K_w) B^{2d_w+2}$ and intuitively assume that FMZVs ``randomly'' behaves in $\Fp$, a non-zero $Q_B$-linear combination of values of $\zeta^{\cA}$ can accidentally give an $\Fp$-linear combination of the corresponding mod $p$ harmonic sums which coincides with $0 \in \Fp$.

\vs
One naive solution is to take a $p \in \bP$ significantly greater than $(\# K_w) B^{2d_w+2}$ in order to avoid the accidental vanishing. However, since Algorithm \ref{parallel horizontal DP} requires $\Theta(p \# V_T)$ time complexity, it is not practical when we consider non-trivial weights. Instead, we apply Chinese remainder algorithm.

\vs
Choose a finite array $\vec{p} = (p_{\ell})_{\ell=0}^{L-1}$ of distinct prime numbers greater than $B$ such that $N \coloneqq \prod_{\ell=0}^{L-1} p_{\ell}$ is significantly greater than $(\# K_w) B^{2d_w+2}$. Then the accidental vanishing in $\prod_{\ell=0}^{L-1} \F_{p_{\ell}} \cong \Z/N \Z$ occurs only in probability approximately $0$. Therefore, computation of $S$ is reduced to the minimal generating system problem for $(\Z/N \Z,\vec{c})$, which is solved by Algorithm \ref{Dynamic MITM algorithm}.

\vs
Computation of mod $p_{\ell}$ harmonic sums for indices in $K_w$ for any $\ell \in \N_{< L}$ is executed in $\Theta(\sum_{\ell=0}^{L-1} p_{\ell} \# V_T)$ time complexity and $\Theta(\# V_T + L \# K_w)$ space complexity by Algorithm \ref{parallel horizontal DP}. There are two way to solve the minimal generating system problem for $(\Z/N \Z,\vec{c})$: One is the application of Garner's algorithm to obtain corresponding values in $\Z/N \Z$, and the other one is to directly handle arrays of length $L$ whose $(1 + \ell)$-th entry is the mod $p_{\ell}$ harmonic sum of a common index $\vec{k}$ for any $\ell \in \N_{< L}$.

\vs
We explain the latter method because it can avoid bignum arithmetic while the former one cannot. A key $y$ of the associative array $\vec{\Sigma}$ for Algorithm \ref{Dynamic MITM algorithm} in the latter method is a tuple of mod $p_{\ell}$ harmonic sums for multiple (but not necessarily all) $\ell \in \N_{< L}$. The optimal number $L^{\rL}$ of chosen $\ell$'s can be dynamically estimated by computing an expected time complexity of the rest of the process, and hence can be updated if necessary when $\vec{\Sigma}$ is reconstructed in the lines 30 -- 34 in Algorithm \ref{Dynamic MITM algorithm}. Here are pseudocodes of the process:

\begin{figure}[H]
\begin{algorithm}[H]
\caption{Pre-computation of for mod $p_{\ell}$ harmonic sums of weight $w$ for each $\ell \in \N_{< L}$ for a $\vec{p} = (p_{\ell})_{\ell=0}^{L-1}$}
\label{Pre-computation}
\begin{algorithmic}[1]
\Function {ModHarmonicSum}{$\vec{p},w$}
	\State $Z \coloneqq ((z_{\ell,\vec{k}})_{\vec{k} \in K_w})_{\ell=0}^{L-1} \gets ((0)_{\vec{k} \in K_w})_{\ell=0}^{L-1}$ 
	\ForAll {$\ell \in \N_{< L}$}
		\State $(v_{\vec{k}})_{\vec{k} \in V_T} \gets$ \Call{parallel horizontal DP}{$p_{\ell},T$}
		\ForAll {$\vec{k} \in K_w$}
			\State $z_{\ell,\vec{k}} \gets v_{\vec{k}}$
		\EndFor
	\EndFor
	\State \Return $Z$
\EndFunction
\end{algorithmic}
\end{algorithm}
\end{figure}

\begin{figure}[H]
\begin{algorithm}[H]
\caption{Expected time complexity for $D^{\rR}$}
\label{Expected time complexity}
\begin{algorithmic}[1]
\Function {CostL}{$B,D,D^{\rL},D^{\rR},H,h,L^{\rL},N$}
	\State \Return $B^{D^{\rL}} L^{\rL} + (H-h) B^{D^{\rR}+1} \lfloor \frac{B^{D^{\rL}}}{N} \rfloor$
\EndFunction
\end{algorithmic}
\end{algorithm}
\end{figure}

\begin{figure}[H]
\begin{algorithm}[H]
\caption{Chinese remainder algorithm for $\vec{p} = (p_{\ell})_{\ell=0}^{L-1} \in \bP_{> B}$}
\label{Chinese remainder algorithm}
\begin{algorithmic}[1]
\Function {ChineseRemainderAgorithm}{$\vec{p},w,B$}
	\State $T \gets$ the tree of indices in Example \ref{tree of indices} (2) applied to $w$
	\State $Z \coloneqq ((z_{\ell,\vec{k}})_{\vec{k} \in K_w})_{\ell=0}^{L-1} \gets$ \Call{ModHarmonicSum}{$\vec{p},w$}
	\State $(L^{\rL},L^{\rL}) \gets (0,L)$ \Comment{the equality $L^{\rR} = L - L^{\rL}$ holds and $L^{\rR}$ is unused}
	\State $N \gets 1$ \Comment{variable for $\prod_{\ell=0}^{L^{\rL}-1} p_{\ell}$}
	\State $(D^{\rL},D^{\rR}) \gets (0,0)$
	\State $X \gets ()$ \Comment{$\vec{x}_d$ for a $d \in \N$ smaller than the current length of $X$ denotes the $(1+d)$-th entry of $X$}
	\State $\vec{\Sigma} \coloneqq (\Sigma_{y'})_{y'} \gets$ the associative array of empty arrays indexed by $\prod_{\ell=0}^{L^{\rL}-1} \F_{p_{\ell}}$
	\State $\vec{c} = (c_b)_{b=0}^{2B} \gets$ an array enumerating $\Z \cap [-B,B]$ in the ascending order of the absolute values
	\State $H \gets \# K_w$
	\State $K = (\vec{k}_h)_{h=0}^{H-1} \gets$ an enumeration of $K_w$
	\ForAll {$h \in \N_{< H}$}
		\If {$(z_{\ell,\vec{k}_h})_{\ell=0}^{L-1}$ is a zero vector}
			\State \Continue
		\EndIf
		\State $\gen \gets$ False
		\ForAll {$c \in \N \cap [1,B]$}
			\ForAll {$(b^{\rR}_d)_{d=0}^{D^{\rR}-1} \in \N_{\leq 2B}^{D_{\rN}}$}
				\State $y \gets ((\sum_{d=0}^{D^{\rR}-1} c_{b^{\rR}_d} z_{\ell,\vec{k}_{D_{\rL}+d}}) + c z_{\ell,\vec{k}_h})_{\ell=0}^{L^{\rL}-1}$
				\ForAll {$\vec{b}^{\rL} \in \Sigma_{-y}$}
					\State $(b_d)_{d=0}^{D-1} \gets \vec{b}^{\rL} \frown \vec{b}^{\rR}$
					\If {$((\sum_{d=0}^{D-1} c_{b_d} z_{\ell,\vec{x}_d}) + c z_{\ell,\vec{k}_h})_{\ell=L^{\rL}}^{L-1}$ is a zero vector}
						\State $\gen \gets$ True
						\State \Break
					\EndIf
				\EndFor
				\If {$\gen$}
					\State \Break
				\EndIf
			\EndFor
			\If {$\gen$}
				\State \Break
			\EndIf
		\EndFor
\algstore{Chinese remainder algorithm 1}
\end{algorithmic}
\end{algorithm}
\end{figure}
\addtocounter{algorithm}{-1}
\begin{figure}[H]
\begin{algorithm}[H]
\begin{algorithmic}[1]
\algrestore{Chinese remainder algorithm 1}
		\If {$\neg \gen$}:
			\State Append $\vec{k}_h$ to $X$
			\State $D \gets D + 1$
			\State $w_0 \gets$ \Call{CostD}{$B,D^{\rR},H,h$}
			\State $w_1 \gets$ \Call{CostD}{$B,D^{\rR}+1,H,h$}
			\If {$w_0 > B^{D^{\rL}+1} + w_1$} \Comment{optimisation of expected computational complexity}
				\State $D^{\rL} \gets D^{\rL} + 1$
				\If {$L^{\rL} < L$}
					\State $w_0 \gets$ \Call{CostL}{$B,D,D^{\rL},D^{\rR},H,h,L^{\rL},N$}
					\State $w_1 \gets$ \Call{CostL}{$B,D,D^{\rL},D^{\rR},H,h,L^{\rL}+1,Np_{L^{\rL}}$}
					\If {$w_0 > w_1$}
						\State $(L^{\rL},L^{\rR}) \gets (L^{\rL}+1,L^{\rR}-1)$
						\State $N \gets N p_{L^{\rL}-1}$
					\EndIf
				\EndIf
				\State $\vec{\Sigma} \gets$ the associative array of empty arrays indexed by $\prod_{\ell=0}^{L^{\rL}-1} \F_{p_{\ell}}$
				\ForAll {$\vec{b}^{\rL} \coloneqq (b_d)_{d=0}^{D^{\rL}-1} \in \N_{\leq 2B}^{D^{\rL}}$}
					\State $y \gets (\sum_{d=0}^{D^{\rL}-1} c_{b_d} z_{\ell,\vec{k}_d})_{\ell=0}^{L^{\rL}-1}$
					\State Append $\vec{b}^{\rL}$ to $\Sigma_y$
				\EndFor
			\Else
				\State $D^{\rR} \gets D^{\rR} + 1$
			\EndIf
		\EndIf
	\EndFor
	\State \Return $X$
\EndFunction
\end{algorithmic}
\end{algorithm}
\end{figure}

\subsection{Experiment}
\label{Experiment}

The deduced $\Q$-linear basis for FMZVs of weight $10$ over $\vec{c}$ given by numerical analysis by \Call{ChineseRemainderAgorithm}{$\vec{p},10,6000$} for
\be
\vec{p} = (10007,10009,10037,10039,10061,10067,10069,10079,10091,10093,10099)
\ee
is $(\zeta^{\cA}(8,1,1),\zeta^{\cA}(7,2,1),\zeta^{\cA}(6,3,1))$. Here is the list of linear relations between the system and other FMZVs of weight $10$:

\begin{table}[H]
\begin{center}
\caption{Linear relations for FMZVs of weight $10$}
\begin{tabular}{|c|c|}
\hline
$\vec{k}$ & coefficients \\
\hline \hline
$(10)$ & $(0,0,0,1)$ \\
\hline
$(9,1)$ & $(0,0,0,1)$ \\
\hline
$(8,2)$ & $(0,0,0,1)$ \\
\hline
$(7,3)$ & $(0,0,0,1)$ \\
\hline
$(7,1,2)$ & $(8,1,0,2)$ \\
\hline
$(7,1,1,1)$ & $(-6,1,0,4)$ \\
\hline
$(6,4)$ & $(0,0,0,1)$ \\
\hline
$(6,2,2)$ & $(0,7,2,2)$ \\
\hline
$(6,2,1,1)$ & $(2,-7,0,4)$ \\
\hline
$(6,1,3)$ & $(-6,-1,0,1)$ \\
\hline
$(6,1,2,1)$ & $(8,3,1,1)$ \\
\hline
$(6,1,1,2)$ & $(-8,-7,-2,2)$ \\
\hline
$(6,1,1,1,1)$ & $(-80,47,2,64)$ \\
\hline
$(5,5)$ & $(0,0,0,1)$ \\
\hline
$(5,4,1)$ & $(14,7,3,1)$ \\
\hline
$(5,3,2)$ & $(-72,-33,-6,2)$ \\
\hline
$(5,3,1,1)$ & $(46,19,0,4)$ \\
\hline
$(5,2,3)$ & $(24,4,0,1)$ \\
\hline
$(5,2,2,1)$ & $(-34,-11,-3,1)$ \\
\hline
$(5,2,1,2)$ & $(19,12,3,1)$ \\
\hline
$(5,2,1,1,1)$ & $(-232,-305,-34,64)$ \\
\hline
$(5,1,4)$ & $(6,1,0,2)$ \\
\hline
$(5,1,3,1)$ & $(-34,-15,-4,4)$ \\
\hline
$(5,1,2,2)$ & $(66,25,4,4)$ \\
\hline
$(5,1,2,1,1)$ & $(752,293,46,64)$ \\
\hline
$(5,1,1,3)$ & $(-42,-7,0,4)$ \\
\hline
$(5,1,1,2,1)$ & $(62,19,13,16)$ \\
\hline
$(5,1,1,1,2)$ & $(-464,-301,-70,64)$ \\
\hline
$(5,1,1,1,1,1)$ & $(-80,47,2,64)$ \\
\hline
$(4,6)$ & $(0,0,0,1)$ \\
\hline
$(4,5,1)$ & $(-34,-15,-6,2)$ \\
\hline
$(4,4,2)$ & $(48,15,2,2)$ \\
\hline
$(4,4,1,1)$ & $(-3,0,1,1)$ \\
\hline
$(4,3,3)$ & $(36,20,4,1)$ \\
\hline
$(4,3,2,1)$ & $(62,17,8,4)$ \\
\hline
$(4,3,1,2)$ & $(-57,-24,-5,1)$ \\
\hline
$(4,3,1,1,1)$ & $(1360,693,110,64)$ \\
\hline
$(4,2,4)$ & $(-48,-15,-2,1)$ \\
\hline
\end{tabular}
\end{center}
\end{table}

\begin{table}[H]
\begin{center}
\caption{Linear relations for FMZVs of weight $10$}
\begin{tabular}{|c|c|}
\hline
$\vec{k}$ & coefficients \\
\hline \hline
$(4,2,3,1)$ & $(130,45,12,4)$ \\
\hline
$(4,2,2,2)$ & $(12,2,0,1)$ \\
\hline
$(4,2,2,1,1)$ & $(-2152,-667,-118,64)$ \\
\hline
$(4,2,1,3)$ & $(9,-2,-1,1)$ \\
\hline
$(4,2,1,2,1)$ & $(-88,-13,-28,16)$ \\
\hline
$(4,2,1,1,2)$ & $(1432,813,186,64)$ \\
\hline
$(4,2,1,1,1,1)$ & $(-232,-305,-34,64)$ \\
\hline
$(4,1,5)$ & $(6,1,0,2)$ \\
\hline
$(4,1,4,1)$ & $(-9,-5,-1,1)$ \\
\hline
$(4,1,3,2)$ & $(90,57,12,4)$ \\
\hline
$(4,1,3,1,1)$ & $(-704,-365,-70,64)$ \\
\hline
$(4,1,2,3)$ & $(-198,-89,-16,4)$ \\
\hline
$(4,1,2,2,1)$ & $(860,401,80,32)$ \\
\hline
$(4,1,2,1,2)$ & $(-1504,-787,-130,64)$ \\
\hline
$(4,1,2,1,1,1)$ & $(208,53,-18,64)$ \\
\hline
$(4,1,1,4)$ & $(30,12,2,1)$ \\
\hline
$(4,1,1,3,1)$ & $(-468,-169,-48,32)$ \\
\hline
$(4,1,1,2,2)$ & $(240,61,6,64)$ \\
\hline
$(4,1,1,2,1,1)$ & $(262,127,29,16)$ \\
\hline
$(4,1,1,1,3)$ & $(240,159,34,32)$ \\
\hline
$(4,1,1,1,2,1)$ & $(-952,-579,-78,64)$ \\
\hline
$(4,1,1,1,1,2)$ & $(50,37,3,16)$ \\
\hline
$(4,1,1,1,1,1,1)$ & $(-6,1,0,4)$ \\
\hline
$(3,7)$ & $(0,0,0,1)$ \\
\hline
$(3,6,1)$ & $(6,1,1,1)$ \\
\hline
$(3,5,2)$ & $(24,25,6,2)$ \\
\hline
$(3,5,1,1)$ & $(-34,-27,-8,4)$ \\
\hline
$(3,4,3)$ & $(-72,-40,-8,1)$ \\
\hline
$(3,4,2,1)$ & $(8,10,1,1)$ \\
\hline
$(3,4,1,2)$ & $(106,47,10,2)$ \\
\hline
$(3,4,1,1,1)$ & $(-314,-193,-35,16)$ \\
\hline
$(3,3,4)$ & $(36,20,4,1)$ \\
\hline
$(3,3,3,1)$ & $(-42,-21,-4,2)$ \\
\hline
$(3,3,2,2)$ & $(-24,-18,-4,1)$ \\
\hline
$(3,3,2,1,1)$ & $(144,119,20,16)$ \\
\hline
$(3,3,1,3)$ & $(42,21,4,2)$ \\
\hline
$(3,3,1,2,1)$ & $(-804,-399,-48,32)$ \\
\hline
$(3,3,1,1,2)$ & $(192,317,62,64)$ \\
\hline
\end{tabular}
\end{center}
\end{table}

\begin{table}[H]
\begin{center}
\caption{Linear relations for FMZVs of weight $10$}
\begin{tabular}{|c|c|}
\hline
$\vec{k}$ & coefficients \\
\hline \hline
$(3,3,1,1,1,1)$ & $(752,293,46,64)$ \\
\hline
$(3,2,5)$ & $(24,4,0,1)$ \\
\hline
$(3,2,4,1)$ & $(-4,4,0,1)$ \\
\hline
$(3,2,3,2)$ & $(-36,-20,-4,1)$ \\
\hline
$(3,2,3,1,1)$ & $(1596,673,144,32)$ \\
\hline
$(3,2,2,3)$ & $(48,36,8,1)$ \\
\hline
$(3,2,2,2,1)$ & $(-336,-447,-66,32)$ \\
\hline
$(3,2,2,1,2)$ & $(-984,-873,-258,64)$ \\
\hline
$(3,2,2,1,1,1)$ & $(152,613,138,64)$ \\
\hline
$(3,2,1,4)$ & $(-198,-89,-16,4)$ \\
\hline
$(3,2,1,3,1)$ & $(1800,783,222,64)$ \\
\hline
$(3,2,1,2,2)$ & $(648,507,114,32)$ \\
\hline
$(3,2,1,2,1,1)$ & $(-1008,-521,-108,16)$ \\
\hline
$(3,2,1,1,3)$ & $(-2208,-1369,-286,64)$ \\
\hline
$(3,2,1,1,2,1)$ & $(2576,1549,206,64)$ \\
\hline
$(3,2,1,1,1,2)$ & $(296,95,38,32)$ \\
\hline
$(3,2,1,1,1,1,1)$ & $(2,-7,0,4)$ \\
\hline
$(3,1,6)$ & $(-6,-1,0,1)$ \\
\hline
$(3,1,5,1)$ & $(34,15,4,4)$ \\
\hline
$(3,1,4,2)$ & $(-30,-19,-4,2)$ \\
\hline
$(3,1,4,1,1)$ & $(-656,-231,-58,64)$ \\
\hline
$(3,1,3,3)$ & $(42,21,4,2)$ \\
\hline
$(3,1,3,2,1)$ & $(280,121,34,64)$ \\
\hline
$(3,1,3,1,2)$ & $(80,23,2,64)$ \\
\hline
$(3,1,3,1,1,1)$ & $(-128,-45,-6,64)$ \\
\hline
$(3,1,2,4)$ & $(9,-2,-1,1)$ \\
\hline
$(3,1,2,3,1)$ & $(-1800,-783,-222,64)$ \\
\hline
$(3,1,2,2,2)$ & $(144,465,126,64)$ \\
\hline
$(3,1,2,2,1,1)$ & $(1020,353,80,32)$ \\
\hline
$(3,1,2,1,3)$ & $(1296,657,126,32)$ \\
\hline
$(3,1,2,1,2,1)$ & $(-312,-141,30,64)$ \\
\hline
$(3,1,2,1,1,2)$ & $(-1460,-753,-152,32)$ \\
\hline
$(3,1,2,1,1,1,1)$ & $(46,19,0,4)$ \\
\hline
$(3,1,1,5)$ & $(-42,-7,0,4)$ \\
\hline
$(3,1,1,4,1)$ & $(210,70,21,16)$ \\
\hline
$(3,1,1,3,2)$ & $(30,19,4,8)$ \\
\hline
$(3,1,1,3,1,1)$ & $(-944,-391,-90,64)$ \\
\hline
$(3,1,1,2,3)$ & $(-2208,-1369,-286,64)$ \\
\hline
\end{tabular}
\end{center}
\end{table}

\begin{table}[H]
\begin{center}
\caption{Linear relations for FMZVs of weight $10$}
\begin{tabular}{|c|c|}
\hline
$\vec{k}$ & coefficients \\
\hline \hline
$(3,1,1,2,2,1)$ & $(-200,5,-22,32)$ \\
\hline
$(3,1,1,2,1,2)$ & $(2576,1363,290,64)$ \\
\hline
$(3,1,1,2,1,1,1)$ & $(-12,-5,0,1)$ \\
\hline
$(3,1,1,1,4)$ & $(240,159,34,32)$ \\
\hline
$(3,1,1,1,3,1)$ & $(24,11,2,32)$ \\
\hline
$(3,1,1,1,2,2)$ & $(-576,-509,-118,64)$ \\
\hline
$(3,1,1,1,2,1,1)$ & $(38,13,0,4)$ \\
\hline
$(3,1,1,1,1,3)$ & $(96,93,22,32)$ \\
\hline
$(3,1,1,1,1,2,1)$ & $(-22,-11,0,4)$ \\
\hline
$(3,1,1,1,1,1,2)$ & $(2,3,0,4)$ \\
\hline
$(3,1,1,1,1,1,1,1)$ & $(-1,0,0,1)$ \\
\hline
$(2,8)$ & $(0,0,0,1)$ \\
\hline
$(2,7,1)$ & $(-8,1,0,2)$ \\
\hline
$(2,6,2)$ & $(0,-7,-2,1)$ \\
\hline
$(2,6,1,1)$ & $(10,13,4,4)$ \\
\hline
$(2,5,3)$ & $(24,25,6,2)$ \\
\hline
$(2,5,2,1)$ & $(11,0,0,1)$ \\
\hline
$(2,5,1,2)$ & $(-23,-9,-2,1)$ \\
\hline
$(2,5,1,1,1)$ & $(272,269,70,64)$ \\
\hline
$(2,4,4)$ & $(48,15,2,2)$ \\
\hline
$(2,4,3,1)$ & $(-10,-4,0,1)$ \\
\hline
$(2,4,2,2)$ & $(-12,-2,0,1)$ \\
\hline
$(2,4,2,1,1)$ & $(728,93,-30,64)$ \\
\hline
$(2,4,1,3)$ & $(-30,-19,-4,2)$ \\
\hline
$(2,4,1,2,1)$ & $(656,515,50,64)$ \\
\hline
$(2,4,1,1,2)$ & $(34,-17,-1,16)$ \\
\hline
$(2,4,1,1,1,1)$ & $(-98,-45,-3,16)$ \\
\hline
$(2,3,5)$ & $(-72,-33,-6,2)$ \\
\hline
$(2,3,4,1)$ & $(-10,-11,-4,2)$ \\
\hline
$(2,3,3,2)$ & $(72,40,8,1)$ \\
\hline
$(2,3,3,1,1)$ & $(-2544,-1213,-206,64)$ \\
\hline
$(2,3,2,3)$ & $(-36,-20,-4,1)$ \\
\hline
$(2,3,2,2,1)$ & $(1848,1059,270,64)$ \\
\hline
$(2,3,2,1,2)$ & $(120,75,30,32)$ \\
\hline
$(2,3,2,1,1,1)$ & $(-184,-141,-44,16)$ \\
\hline
$(2,3,1,4)$ & $(90,57,12,4)$ \\
\hline
$(2,3,1,3,1)$ & $(-90,-36,-9,16)$ \\
\hline
$(2,3,1,2,2)$ & $(-864,-921,-222,64)$ \\
\hline
\end{tabular}
\end{center}
\end{table}

\begin{table}[H]
\begin{center}
\caption{Linear relations for FMZVs of weight $10$}
\begin{tabular}{|c|c|}
\hline
$\vec{k}$ & coefficients \\
\hline \hline
$(2,3,1,2,1,1)$ & $(860,625,144,32)$ \\
\hline
$(2,3,1,1,3)$ & $(30,19,4,8)$ \\
\hline
$(2,3,1,1,2,1)$ & $(-1592,-1101,-226,64)$ \\
\hline
$(2,3,1,1,1,2)$ & $(188,231,40,32)$ \\
\hline
$(2,3,1,1,1,1,1)$ & $(-2,-1,0,1)$ \\
\hline
$(2,2,6)$ & $(0,7,2,2)$ \\
\hline
$(2,2,5,1)$ & $(42,7,4,4)$ \\
\hline
$(2,2,4,2)$ & $(-12,-2,0,1)$ \\
\hline
$(2,2,4,1,1)$ & $(28,29,-4,32)$ \\
\hline
$(2,2,3,3)$ & $(-24,-18,-4,1)$ \\
\hline
$(2,2,3,2,1)$ & $(-752,-227,-86,32)$ \\
\hline
$(2,2,3,1,2)$ & $(140,83,32,32)$ \\
\hline
$(2,2,3,1,1,1)$ & $(316,49,16,32)$ \\
\hline
$(2,2,2,4)$ & $(12,2,0,1)$ \\
\hline
$(2,2,2,3,1)$ & $(54,72,9,16)$ \\
\hline
$(2,2,2,2,2)$ & $(0,0,0,1)$ \\
\hline
$(2,2,2,2,1,1)$ & $(112,1,-2,32)$ \\
\hline
$(2,2,2,1,3)$ & $(144,465,126,64)$ \\
\hline
$(2,2,2,1,2,1)$ & $(-528,-429,-6,64)$ \\
\hline
$(2,2,2,1,1,2)$ & $(80,23,2,8)$ \\
\hline
$(2,2,2,1,1,1,1)$ & $(-4,1,0,1)$ \\
\hline
$(2,2,1,5)$ & $(66,25,4,4)$ \\
\hline
$(2,2,1,4,1)$ & $(-792,-153,-42,64)$ \\
\hline
$(2,2,1,3,2)$ & $(-864,-921,-222,64)$ \\
\hline
$(2,2,1,3,1,1)$ & $(1144,601,130,64)$ \\
\hline
$(2,2,1,2,3)$ & $(648,507,114,32)$ \\
\hline
$(2,2,1,2,2,1)$ & $(492,69,24,32)$ \\
\hline
$(2,2,1,2,1,2)$ & $(312,141,-30,64)$ \\
\hline
$(2,2,1,2,1,1,1)$ & $(-30,-3,0,4)$ \\
\hline
$(2,2,1,1,4)$ & $(240,61,6,64)$ \\
\hline
$(2,2,1,1,3,1)$ & $(-160,-169,2,64)$ \\
\hline
$(2,2,1,1,2,2)$ & $(-432,-93,-6,16)$ \\
\hline
$(2,2,1,1,2,1,1)$ & $(4,0,0,1)$ \\
\hline
$(2,2,1,1,1,3)$ & $(-576,-509,-118,64)$ \\
\hline
$(2,2,1,1,1,2,1)$ & $(5,6,0,1)$ \\
\hline
$(2,2,1,1,1,1,2)$ & $(4,-1,0,1)$ \\
\hline
$(2,2,1,1,1,1,1,1)$ & $(0,-1,0,1)$ \\
\hline
$(2,1,7)$ & $(8,1,0,2)$ \\
\hline
\end{tabular}
\end{center}
\end{table}

\begin{table}[H]
\begin{center}
\caption{Linear relations for FMZVs of weight $10$}
\begin{tabular}{|c|c|}
\hline
$\vec{k}$ & coefficients \\
\hline \hline
$(2,1,6,1)$ & $(2,1,0,1)$ \\
\hline
$(2,1,5,2)$ & $(-23,-9,-2,1)$ \\
\hline
$(2,1,5,1,1)$ & $(344,139,34,32)$ \\
\hline
$(2,1,4,3)$ & $(106,47,10,2)$ \\
\hline
$(2,1,4,2,1)$ & $(-1440,-681,-126,64)$ \\
\hline
$(2,1,4,1,2)$ & $(160,55,10,32)$ \\
\hline
$(2,1,4,1,1,1)$ & $(72,47,6,32)$ \\
\hline
$(2,1,3,4)$ & $(-57,-24,-5,1)$ \\
\hline
$(2,1,3,3,1)$ & $(1592,653,98,64)$ \\
\hline
$(2,1,3,2,2)$ & $(140,83,32,32)$ \\
\hline
$(2,1,3,2,1,1)$ & $(-960,-321,-54,64)$ \\
\hline
$(2,1,3,1,3)$ & $(80,23,2,64)$ \\
\hline
$(2,1,3,1,2,1)$ & $(270,108,27,16)$ \\
\hline
$(2,1,3,1,1,2)$ & $(-400,-169,-46,32)$ \\
\hline
$(2,1,3,1,1,1,1)$ & $(6,1,0,4)$ \\
\hline
$(2,1,2,5)$ & $(19,12,3,1)$ \\
\hline
$(2,1,2,4,1)$ & $(80,-19,62,64)$ \\
\hline
$(2,1,2,3,2)$ & $(120,75,30,32)$ \\
\hline
$(2,1,2,3,1,1)$ & $(-720,-351,-114,64)$ \\
\hline
$(2,1,2,2,3)$ & $(-984,-873,-258,64)$ \\
\hline
$(2,1,2,2,2,1)$ & $(-306,-72,-27,16)$ \\
\hline
$(2,1,2,2,1,2)$ & $(120,75,30,16)$ \\
\hline
$(2,1,2,2,1,1,1)$ & $(10,-3,0,4)$ \\
\hline
$(2,1,2,1,4)$ & $(-1504,-787,-130,64)$ \\
\hline
$(2,1,2,1,3,1)$ & $(672,285,6,64)$ \\
\hline
$(2,1,2,1,2,2)$ & $(312,141,-30,64)$ \\
\hline
$(2,1,2,1,2,1,1)$ & $(-2,-5,0,2)$ \\
\hline
$(2,1,2,1,1,3)$ & $(2576,1363,290,64)$ \\
\hline
$(2,1,2,1,1,2,1)$ & $(-10,-4,0,1)$ \\
\hline
$(2,1,2,1,1,1,2)$ & $(-50,-17,0,4)$ \\
\hline
$(2,1,2,1,1,1,1,1)$ & $(10,4,0,1)$ \\
\hline
$(2,1,1,6)$ & $(-8,-7,-2,2)$ \\
\hline
$(2,1,1,5,1)$ & $(-200,-5,-34,64)$ \\
\hline
$(2,1,1,4,2)$ & $(34,-17,-1,16)$ \\
\hline
$(2,1,1,4,1,1)$ & $(12,9,2,8)$ \\
\hline
$(2,1,1,3,3)$ & $(192,317,62,64)$ \\
\hline
$(2,1,1,3,2,1)$ & $(1320,423,150,64)$ \\
\hline
$(2,1,1,3,1,2)$ & $(-400,-169,-46,32)$ \\
\hline
\end{tabular}
\end{center}
\end{table}

\begin{table}[H]
\begin{center}
\caption{Linear relations for FMZVs of weight $10$}
\begin{tabular}{|c|c|}
\hline
$\vec{k}$ & coefficients \\
\hline \hline
$(2,1,1,3,1,1,1)$ & $(-1,1,0,1)$ \\
\hline
$(2,1,1,2,4)$ & $(1432,813,186,64)$ \\
\hline
$(2,1,1,2,3,1)$ & $(-1160,-515,-62,64)$ \\
\hline
$(2,1,1,2,2,2)$ & $(80,23,2,8)$ \\
\hline
$(2,1,1,2,2,1,1)$ & $(0,0,0,1)$ \\
\hline
$(2,1,1,2,1,3)$ & $(-1460,-753,-152,32)$ \\
\hline
$(2,1,1,2,1,2,1)$ & $(22,13,0,2)$ \\
\hline
$(2,1,1,2,1,1,2)$ & $(16,4,0,1)$ \\
\hline
$(2,1,1,2,1,1,1,1)$ & $(-16,-5,0,1)$ \\
\hline
$(2,1,1,1,5)$ & $(-464,-301,-70,64)$ \\
\hline
$(2,1,1,1,4,1)$ & $(-32,-117,-30,64)$ \\
\hline
$(2,1,1,1,3,2)$ & $(188,231,40,32)$ \\
\hline
$(2,1,1,1,3,1,1)$ & $(-2,-5,0,4)$ \\
\hline
$(2,1,1,1,2,3)$ & $(296,95,38,32)$ \\
\hline
$(2,1,1,1,2,2,1)$ & $(-10,3,0,4)$ \\
\hline
$(2,1,1,1,2,1,2)$ & $(-50,-17,0,4)$ \\
\hline
$(2,1,1,1,2,1,1,1)$ & $(30,5,0,2)$ \\
\hline
$(2,1,1,1,1,4)$ & $(50,37,3,16)$ \\
\hline
$(2,1,1,1,1,3,1)$ & $(2,1,0,1)$ \\
\hline
$(2,1,1,1,1,2,2)$ & $(4,-1,0,1)$ \\
\hline
$(2,1,1,1,1,2,1,1)$ & $(-8,1,0,1)$ \\
\hline
$(2,1,1,1,1,1,3)$ & $(2,3,0,4)$ \\
\hline
$(2,1,1,1,1,1,2,1)$ & $(0,-7,0,2)$ \\
\hline
$(2,1,1,1,1,1,1,2)$ & $(0,2,0,1)$ \\
\hline
$(2,1,1,1,1,1,1,1,1)$ & $(0,0,0,1)$ \\
\hline
$(1,9)$ & $(0,0,0,1)$ \\
\hline
$(1,8,1)$ & $(2,0,0,1)$ \\
\hline
$(1,7,2)$ & $(-8,1,0,2)$ \\
\hline
$(1,7,1,1)$ & $(6,-1,0,4)$ \\
\hline
$(1,6,3)$ & $(6,1,1,1)$ \\
\hline
$(1,6,2,1)$ & $(-9,-2,-1,1)$ \\
\hline
$(1,6,1,2)$ & $(2,1,0,1)$ \\
\hline
$(1,6,1,1,1)$ & $(32,-15,-2,16)$ \\
\hline
$(1,5,4)$ & $(-34,-15,-6,2)$ \\
\hline
$(1,5,3,1)$ & $(30,15,8,4)$ \\
\hline
$(1,5,2,2)$ & $(42,7,4,4)$ \\
\hline
$(1,5,2,1,1)$ & $(68,183,28,32)$ \\
\hline
$(1,5,1,3)$ & $(34,15,4,4)$ \\
\hline
\end{tabular}
\end{center}
\end{table}

\begin{table}[H]
\begin{center}
\caption{Linear relations for FMZVs of weight $10$}
\begin{tabular}{|c|c|}
\hline
$\vec{k}$ & coefficients \\
\hline \hline
$(1,5,1,2,1)$ & $(-424,-193,-34,32)$ \\
\hline
$(1,5,1,1,2)$ & $(-200,-5,-34,64)$ \\
\hline
$(1,5,1,1,1,1)$ & $(176,-45,-6,64)$ \\
\hline
$(1,4,5)$ & $(14,7,3,1)$ \\
\hline
$(1,4,4,1)$ & $(-6,-2,-2,1)$ \\
\hline
$(1,4,3,2)$ & $(-10,-11,-4,2)$ \\
\hline
$(1,4,3,1,1)$ & $(-632,-377,-66,32)$ \\
\hline
$(1,4,2,3)$ & $(-4,4,0,1)$ \\
\hline
$(1,4,2,2,1)$ & $(318,42,15,16)$ \\
\hline
$(1,4,2,1,2)$ & $(80,-19,62,64)$ \\
\hline
$(1,4,2,1,1,1)$ & $(664,685,122,64)$ \\
\hline
$(1,4,1,4)$ & $(-9,-5,-1,1)$ \\
\hline
$(1,4,1,3,1)$ & $(440,167,38,32)$ \\
\hline
$(1,4,1,2,2)$ & $(-792,-153,-42,64)$ \\
\hline
$(1,4,1,2,1,1)$ & $(-832,-451,-66,64)$ \\
\hline
$(1,4,1,1,3)$ & $(210,70,21,16)$ \\
\hline
$(1,4,1,1,2,1)$ & $(-40,43,-2,16)$ \\
\hline
$(1,4,1,1,1,2)$ & $(-32,-117,-30,64)$ \\
\hline
$(1,4,1,1,1,1,1)$ & $(12,1,0,2)$ \\
\hline
$(1,3,6)$ & $(0,0,-1,1)$ \\
\hline
$(1,3,5,1)$ & $(30,15,8,4)$ \\
\hline
$(1,3,4,2)$ & $(-10,-4,0,1)$ \\
\hline
$(1,3,4,1,1)$ & $(1256,835,158,64)$ \\
\hline
$(1,3,3,3)$ & $(-42,-21,-4,2)$ \\
\hline
$(1,3,3,2,1)$ & $(1552,637,94,64)$ \\
\hline
$(1,3,3,1,2)$ & $(1592,653,98,64)$ \\
\hline
$(1,3,3,1,1,1)$ & $(-2144,-1043,-194,64)$ \\
\hline
$(1,3,2,4)$ & $(130,45,12,4)$ \\
\hline
$(1,3,2,3,1)$ & $(-1140,-429,-96,16)$ \\
\hline
$(1,3,2,2,2)$ & $(54,72,9,16)$ \\
\hline
$(1,3,2,2,1,1)$ & $(1576,151,-2,64)$ \\
\hline
$(1,3,2,1,3)$ & $(-1800,-783,-222,64)$ \\
\hline
$(1,3,2,1,2,1)$ & $(1380,579,156,32)$ \\
\hline
$(1,3,2,1,1,2)$ & $(-1160,-515,-62,64)$ \\
\hline
$(1,3,2,1,1,1,1)$ & $(-11,0,0,1)$ \\
\hline
$(1,3,1,5)$ & $(-34,-15,-4,4)$ \\
\hline
$(1,3,1,4,1)$ & $(440,167,38,32)$ \\
\hline
$(1,3,1,3,2)$ & $(-90,-36,-9,16)$ \\
\hline
\end{tabular}
\end{center}
\end{table}

\begin{table}[H]
\begin{center}
\caption{Linear relations for FMZVs of weight $10$}
\begin{tabular}{|c|c|}
\hline
$\vec{k}$ & coefficients \\
\hline \hline
$(1,3,1,3,1,1)$ & $(80,23,2,64)$ \\
\hline
$(1,3,1,2,3)$ & $(1800,783,222,64)$ \\
\hline
$(1,3,1,2,2,1)$ & $(-1320,-501,-114,32)$ \\
\hline
$(1,3,1,2,1,2)$ & $(672,285,6,64)$ \\
\hline
$(1,3,1,2,1,1,1)$ & $(-7,-4,0,1)$ \\
\hline
$(1,3,1,1,4)$ & $(-468,-169,-48,32)$ \\
\hline
$(1,3,1,1,3,1)$ & $(800,311,74,32)$ \\
\hline
$(1,3,1,1,2,2)$ & $(-160,-169,2,64)$ \\
\hline
$(1,3,1,1,2,1,1)$ & $(6,7,0,2)$ \\
\hline
$(1,3,1,1,1,3)$ & $(24,11,2,32)$ \\
\hline
$(1,3,1,1,1,2,1)$ & $(-3,-1,0,1)$ \\
\hline
$(1,3,1,1,1,1,2)$ & $(2,1,0,1)$ \\
\hline
$(1,3,1,1,1,1,1,1)$ & $(8,1,0,2)$ \\
\hline
$(1,2,7)$ & $(0,-1,0,1)$ \\
\hline
$(1,2,6,1)$ & $(-9,-2,-1,1)$ \\
\hline
$(1,2,5,2)$ & $(11,0,0,1)$ \\
\hline
$(1,2,5,1,1)$ & $(-352,-285,-86,64)$ \\
\hline
$(1,2,4,3)$ & $(8,10,1,1)$ \\
\hline
$(1,2,4,2,1)$ & $(-120,-69,18,16)$ \\
\hline
$(1,2,4,1,2)$ & $(-1440,-681,-126,64)$ \\
\hline
$(1,2,4,1,1,1)$ & $(888,589,90,64)$ \\
\hline
$(1,2,3,4)$ & $(62,17,8,4)$ \\
\hline
$(1,2,3,3,1)$ & $(1552,637,94,64)$ \\
\hline
$(1,2,3,2,2)$ & $(-752,-227,-86,32)$ \\
\hline
$(1,2,3,2,1,1)$ & $(300,99,60,32)$ \\
\hline
$(1,2,3,1,3)$ & $(280,121,34,64)$ \\
\hline
$(1,2,3,1,2,1)$ & $(-1008,-441,-198,64)$ \\
\hline
$(1,2,3,1,1,2)$ & $(1320,423,150,64)$ \\
\hline
$(1,2,3,1,1,1,1)$ & $(26,9,0,4)$ \\
\hline
$(1,2,2,5)$ & $(-34,-11,-3,1)$ \\
\hline
$(1,2,2,4,1)$ & $(318,42,15,16)$ \\
\hline
$(1,2,2,3,2)$ & $(1848,1059,270,64)$ \\
\hline
$(1,2,2,3,1,1)$ & $(-1800,-447,-78,64)$ \\
\hline
$(1,2,2,2,3)$ & $(-336,-447,-66,32)$ \\
\hline
$(1,2,2,2,2,1)$ & $(54,72,9,8)$ \\
\hline
$(1,2,2,2,1,2)$ & $(-306,-72,-27,16)$ \\
\hline
$(1,2,2,2,1,1,1)$ & $(15,-1,0,1)$ \\
\hline
$(1,2,2,1,4)$ & $(860,401,80,32)$ \\
\hline
\end{tabular}
\end{center}
\end{table}

\begin{table}[H]
\begin{center}
\caption{Linear relations for FMZVs of weight $10$}
\begin{tabular}{|c|c|}
\hline
$\vec{k}$ & coefficients \\
\hline \hline
$(1,2,2,1,3,1)$ & $(-1320,-501,-114,32)$ \\
\hline
$(1,2,2,1,2,2)$ & $(492,69,24,32)$ \\
\hline
$(1,2,2,1,2,1,1)$ & $(10,4,0,1)$ \\
\hline
$(1,2,2,1,1,3)$ & $(-200,5,-22,32)$ \\
\hline
$(1,2,2,1,1,2,1)$ & $(-6,-8,0,1)$ \\
\hline
$(1,2,2,1,1,1,2)$ & $(-10,3,0,4)$ \\
\hline
$(1,2,2,1,1,1,1,1)$ & $(-20,-1,0,2)$ \\
\hline
$(1,2,1,6)$ & $(8,3,1,1)$ \\
\hline
$(1,2,1,5,1)$ & $(-424,-193,-34,32)$ \\
\hline
$(1,2,1,4,2)$ & $(656,515,50,64)$ \\
\hline
$(1,2,1,4,1,1)$ & $(102,10,9,16)$ \\
\hline
$(1,2,1,3,3)$ & $(-804,-399,-48,32)$ \\
\hline
$(1,2,1,3,2,1)$ & $(-1008,-441,-198,64)$ \\
\hline
$(1,2,1,3,1,2)$ & $(270,108,27,16)$ \\
\hline
$(1,2,1,3,1,1,1)$ & $(-14,-7,0,4)$ \\
\hline
$(1,2,1,2,4)$ & $(-88,-13,-28,16)$ \\
\hline
$(1,2,1,2,3,1)$ & $(1380,579,156,32)$ \\
\hline
$(1,2,1,2,2,2)$ & $(-528,-429,-6,64)$ \\
\hline
$(1,2,1,2,2,1,1)$ & $(-4,4,0,1)$ \\
\hline
$(1,2,1,2,1,3)$ & $(-312,-141,30,64)$ \\
\hline
$(1,2,1,2,1,2,1)$ & $(0,0,0,1)$ \\
\hline
$(1,2,1,2,1,1,2)$ & $(22,13,0,2)$ \\
\hline
$(1,2,1,2,1,1,1,1)$ & $(-12,-9,0,2)$ \\
\hline
$(1,2,1,1,5)$ & $(62,19,13,16)$ \\
\hline
$(1,2,1,1,4,1)$ & $(-40,43,-2,16)$ \\
\hline
$(1,2,1,1,3,2)$ & $(-1592,-1101,-226,64)$ \\
\hline
$(1,2,1,1,3,1,1)$ & $(6,1,0,2)$ \\
\hline
$(1,2,1,1,2,3)$ & $(2576,1549,206,64)$ \\
\hline
$(1,2,1,1,2,2,1)$ & $(-6,-8,0,1)$ \\
\hline
$(1,2,1,1,2,1,2)$ & $(-10,-4,0,1)$ \\
\hline
$(1,2,1,1,2,1,1,1)$ & $(20,15,0,2)$ \\
\hline
$(1,2,1,1,1,4)$ & $(-952,-579,-78,64)$ \\
\hline
$(1,2,1,1,1,3,1)$ & $(-3,-1,0,1)$ \\
\hline
$(1,2,1,1,1,2,2)$ & $(5,6,0,1)$ \\
\hline
$(1,2,1,1,1,2,1,1)$ & $(-20,-15,0,2)$ \\
\hline
$(1,2,1,1,1,1,3)$ & $(-22,-11,0,4)$ \\
\hline
$(1,2,1,1,1,1,2,1)$ & $(12,9,0,1)$ \\
\hline
$(1,2,1,1,1,1,1,2)$ & $(0,-7,0,2)$ \\
\hline
\end{tabular}
\end{center}
\end{table}

\begin{table}[H]
\begin{center}
\caption{Linear relations for FMZVs of weight $10$}
\begin{tabular}{|c|c|}
\hline
$\vec{k}$ & coefficients \\
\hline \hline
$(1,2,1,1,1,1,1,1,1)$ & $(0,0,0,1)$ \\
\hline
$(1,1,8)$ & $(-1,0,0,1)$ \\
\hline
$(1,1,7,1)$ & $(6,-1,0,4)$ \\
\hline
$(1,1,6,2)$ & $(10,13,4,4)$ \\
\hline
$(1,1,6,1,1)$ & $(-48,13,6,32)$ \\
\hline
$(1,1,5,3)$ & $(-34,-27,-8,4)$ \\
\hline
$(1,1,5,2,1)$ & $(-352,-285,-86,64)$ \\
\hline
$(1,1,5,1,2)$ & $(344,139,34,32)$ \\
\hline
$(1,1,5,1,1,1)$ & $(-48,-1,2,32)$ \\
\hline
$(1,1,4,4)$ & $(-3,0,1,1)$ \\
\hline
$(1,1,4,3,1)$ & $(1256,835,158,64)$ \\
\hline
$(1,1,4,2,2)$ & $(28,29,-4,32)$ \\
\hline
$(1,1,4,2,1,1)$ & $(-1056,-729,-158,64)$ \\
\hline
$(1,1,4,1,3)$ & $(-656,-231,-58,64)$ \\
\hline
$(1,1,4,1,2,1)$ & $(102,10,9,16)$ \\
\hline
$(1,1,4,1,1,2)$ & $(12,9,2,8)$ \\
\hline
$(1,1,4,1,1,1,1)$ & $(-42,-7,0,4)$ \\
\hline
$(1,1,3,5)$ & $(46,19,0,4)$ \\
\hline
$(1,1,3,4,1)$ & $(-632,-377,-66,32)$ \\
\hline
$(1,1,3,3,2)$ & $(-2544,-1213,-206,64)$ \\
\hline
$(1,1,3,3,1,1)$ & $(1440,737,142,32)$ \\
\hline
$(1,1,3,2,3)$ & $(1596,673,144,32)$ \\
\hline
$(1,1,3,2,2,1)$ & $(-1800,-447,-78,64)$ \\
\hline
$(1,1,3,2,1,2)$ & $(-720,-351,-114,64)$ \\
\hline
$(1,1,3,2,1,1,1)$ & $(18,3,0,1)$ \\
\hline
$(1,1,3,1,4)$ & $(-704,-365,-70,64)$ \\
\hline
$(1,1,3,1,3,1)$ & $(80,23,2,64)$ \\
\hline
$(1,1,3,1,2,2)$ & $(1144,601,130,64)$ \\
\hline
$(1,1,3,1,2,1,1)$ & $(6,1,0,2)$ \\
\hline
$(1,1,3,1,1,3)$ & $(-944,-391,-90,64)$ \\
\hline
$(1,1,3,1,1,2,1)$ & $(6,1,0,2)$ \\
\hline
$(1,1,3,1,1,1,2)$ & $(-2,-5,0,4)$ \\
\hline
$(1,1,3,1,1,1,1,1)$ & $(-6,-1,0,1)$ \\
\hline
$(1,1,2,6)$ & $(2,-7,0,4)$ \\
\hline
$(1,1,2,5,1)$ & $(68,183,28,32)$ \\
\hline
$(1,1,2,4,2)$ & $(728,93,-30,64)$ \\
\hline
$(1,1,2,4,1,1)$ & $(-1056,-729,-158,64)$ \\
\hline
$(1,1,2,3,3)$ & $(144,119,20,16)$ \\
\hline
\end{tabular}
\end{center}
\end{table}

\begin{table}[H]
\begin{center}
\caption{Linear relations for FMZVs of weight $10$}
\begin{tabular}{|c|c|}
\hline
$\vec{k}$ & coefficients \\
\hline \hline
$(1,1,2,3,2,1)$ & $(300,99,60,32)$ \\
\hline
$(1,1,2,3,1,2)$ & $(-960,-321,-54,64)$ \\
\hline
$(1,1,2,3,1,1,1)$ & $(-54,-9,0,4)$ \\
\hline
$(1,1,2,2,4)$ & $(-2152,-667,-118,64)$ \\
\hline
$(1,1,2,2,3,1)$ & $(1576,151,-2,64)$ \\
\hline
$(1,1,2,2,2,2)$ & $(112,1,-2,32)$ \\
\hline
$(1,1,2,2,2,1,1)$ & $(-24,-4,0,1)$ \\
\hline
$(1,1,2,2,1,3)$ & $(1020,353,80,32)$ \\
\hline
$(1,1,2,2,1,2,1)$ & $(-4,4,0,1)$ \\
\hline
$(1,1,2,2,1,1,2)$ & $(0,0,0,1)$ \\
\hline
$(1,1,2,2,1,1,1,1)$ & $(24,4,0,1)$ \\
\hline
$(1,1,2,1,5)$ & $(752,293,46,64)$ \\
\hline
$(1,1,2,1,4,1)$ & $(-832,-451,-66,64)$ \\
\hline
$(1,1,2,1,3,2)$ & $(860,625,144,32)$ \\
\hline
$(1,1,2,1,3,1,1)$ & $(6,1,0,2)$ \\
\hline
$(1,1,2,1,2,3)$ & $(-1008,-521,-108,16)$ \\
\hline
$(1,1,2,1,2,2,1)$ & $(10,4,0,1)$ \\
\hline
$(1,1,2,1,2,1,2)$ & $(-2,-5,0,2)$ \\
\hline
$(1,1,2,1,2,1,1,1)$ & $(0,0,0,1)$ \\
\hline
$(1,1,2,1,1,4)$ & $(262,127,29,16)$ \\
\hline
$(1,1,2,1,1,3,1)$ & $(6,7,0,2)$ \\
\hline
$(1,1,2,1,1,2,2)$ & $(4,0,0,1)$ \\
\hline
$(1,1,2,1,1,2,1,1)$ & $(0,0,0,1)$ \\
\hline
$(1,1,2,1,1,1,3)$ & $(38,13,0,4)$ \\
\hline
$(1,1,2,1,1,1,2,1)$ & $(-20,-15,0,2)$ \\
\hline
$(1,1,2,1,1,1,1,2)$ & $(-8,1,0,1)$ \\
\hline
$(1,1,2,1,1,1,1,1,1)$ & $(0,0,0,1)$ \\
\hline
$(1,1,1,7)$ & $(-6,1,0,4)$ \\
\hline
$(1,1,1,6,1)$ & $(32,-15,-2,16)$ \\
\hline
$(1,1,1,5,2)$ & $(272,269,70,64)$ \\
\hline
$(1,1,1,5,1,1)$ & $(-48,-1,2,32)$ \\
\hline
$(1,1,1,4,3)$ & $(-314,-193,-35,16)$ \\
\hline
$(1,1,1,4,2,1)$ & $(888,589,90,64)$ \\
\hline
$(1,1,1,4,1,2)$ & $(72,47,6,32)$ \\
\hline
$(1,1,1,4,1,1,1)$ & $(12,2,0,1)$ \\
\hline
$(1,1,1,3,4)$ & $(1360,693,110,64)$ \\
\hline
$(1,1,1,3,3,1)$ & $(-2144,-1043,-194,64)$ \\
\hline
$(1,1,1,3,2,2)$ & $(316,49,16,32)$ \\
\hline
\end{tabular}
\end{center}
\end{table}

\begin{table}[H]
\begin{center}
\caption{Linear relations for FMZVs of weight $10$}
\begin{tabular}{|c|c|}
\hline
$\vec{k}$ & coefficients \\
\hline \hline
$(1,1,1,3,2,1,1)$ & $(-54,-9,0,4)$ \\
\hline
$(1,1,1,3,1,3)$ & $(-128,-45,-6,64)$ \\
\hline
$(1,1,1,3,1,2,1)$ & $(-14,-7,0,4)$ \\
\hline
$(1,1,1,3,1,1,2)$ & $(-1,1,0,1)$ \\
\hline
$(1,1,1,3,1,1,1,1)$ & $(6,1,0,2)$ \\
\hline
$(1,1,1,2,5)$ & $(-232,-305,-34,64)$ \\
\hline
$(1,1,1,2,4,1)$ & $(664,685,122,64)$ \\
\hline
$(1,1,1,2,3,2)$ & $(-184,-141,-44,16)$ \\
\hline
$(1,1,1,2,3,1,1)$ & $(18,3,0,1)$ \\
\hline
$(1,1,1,2,2,3)$ & $(152,613,138,64)$ \\
\hline
$(1,1,1,2,2,2,1)$ & $(15,-1,0,1)$ \\
\hline
$(1,1,1,2,2,1,2)$ & $(10,-3,0,4)$ \\
\hline
$(1,1,1,2,2,1,1,1)$ & $(-30,-5,0,1)$ \\
\hline
$(1,1,1,2,1,4)$ & $(208,53,-18,64)$ \\
\hline
$(1,1,1,2,1,3,1)$ & $(-7,-4,0,1)$ \\
\hline
$(1,1,1,2,1,2,2)$ & $(-30,-3,0,4)$ \\
\hline
$(1,1,1,2,1,2,1,1)$ & $(0,0,0,1)$ \\
\hline
$(1,1,1,2,1,1,3)$ & $(-12,-5,0,1)$ \\
\hline
$(1,1,1,2,1,1,2,1)$ & $(20,15,0,2)$ \\
\hline
$(1,1,1,2,1,1,1,2)$ & $(30,5,0,2)$ \\
\hline
$(1,1,1,2,1,1,1,1,1)$ & $(0,0,0,1)$ \\
\hline
$(1,1,1,1,6)$ & $(-80,47,2,64)$ \\
\hline
$(1,1,1,1,5,1)$ & $(176,-45,-6,64)$ \\
\hline
$(1,1,1,1,4,2)$ & $(-98,-45,-3,16)$ \\
\hline
$(1,1,1,1,4,1,1)$ & $(-42,-7,0,4)$ \\
\hline
$(1,1,1,1,3,3)$ & $(752,293,46,64)$ \\
\hline
$(1,1,1,1,3,2,1)$ & $(26,9,0,4)$ \\
\hline
$(1,1,1,1,3,1,2)$ & $(6,1,0,4)$ \\
\hline
$(1,1,1,1,3,1,1,1)$ & $(6,1,0,2)$ \\
\hline
$(1,1,1,1,2,4)$ & $(-232,-305,-34,64)$ \\
\hline
$(1,1,1,1,2,3,1)$ & $(-11,0,0,1)$ \\
\hline
$(1,1,1,1,2,2,2)$ & $(-4,1,0,1)$ \\
\hline
$(1,1,1,1,2,2,1,1)$ & $(24,4,0,1)$ \\
\hline
$(1,1,1,1,2,1,3)$ & $(46,19,0,4)$ \\
\hline
$(1,1,1,1,2,1,2,1)$ & $(-12,-9,0,2)$ \\
\hline
$(1,1,1,1,2,1,1,2)$ & $(-16,-5,0,1)$ \\
\hline
$(1,1,1,1,2,1,1,1,1)$ & $(0,0,0,1)$ \\
\hline
$(1,1,1,1,1,5)$ & $(-80,47,2,64)$ \\
\hline
\end{tabular}
\end{center}
\end{table}

\begin{table}[H]
\begin{center}
\caption{Linear relations for FMZVs of weight $10$}
\begin{tabular}{|c|c|}
\hline
$\vec{k}$ & coefficients \\
\hline \hline
$(1,1,1,1,1,4,1)$ & $(12,1,0,2)$ \\
\hline
$(1,1,1,1,1,3,2)$ & $(-2,-1,0,1)$ \\
\hline
$(1,1,1,1,1,3,1,1)$ & $(-6,-1,0,1)$ \\
\hline
$(1,1,1,1,1,2,3)$ & $(2,-7,0,4)$ \\
\hline
$(1,1,1,1,1,2,2,1)$ & $(-20,-1,0,2)$ \\
\hline
$(1,1,1,1,1,2,1,2)$ & $(10,4,0,1)$ \\
\hline
$(1,1,1,1,1,2,1,1,1)$ & $(0,0,0,1)$ \\
\hline
$(1,1,1,1,1,1,4)$ & $(-6,1,0,4)$ \\
\hline
$(1,1,1,1,1,1,3,1)$ & $(8,1,0,2)$ \\
\hline
$(1,1,1,1,1,1,2,2)$ & $(0,-1,0,1)$ \\
\hline
$(1,1,1,1,1,1,2,1,1)$ & $(0,0,0,1)$ \\
\hline
$(1,1,1,1,1,1,1,3)$ & $(-1,0,0,1)$ \\
\hline
$(1,1,1,1,1,1,1,2,1)$ & $(0,0,0,1)$ \\
\hline
$(1,1,1,1,1,1,1,1,2)$ & $(0,0,0,1)$ \\
\hline
$(1,1,1,1,1,1,1,1,1,1)$ & $(0,0,0,1)$ \\
\hline
\end{tabular}
\end{center}
\end{table}

For example, the line with $(6,1,1,1,1)$ and $(-80,47,2,64)$ means that the linear combination
\be
-80 \zeta^{\cA}(8,1,1) + 47 \zeta^{\cA}(7,2,1) + 2 \zeta^{\cA}(6,3,1) + 64 \zeta^{\cA}(6,1,1,1,1)
\ee
is deduced to be $0$ from the numerical analysis. Since
\be
N & = & 10007 \times 10009 \times 10037 \times 10039 \times 10061 \times 10067 \\
& & \times 10069 \times 10079 \times 10091 \times 10093 \times 10099 \\
& = & 106700590455862347842907841856033238416352421 \approx 10^{44}
\ee
is significantly larger than
\be
\# K_{10} \times 6000^3 = 512 \times 6000^3 = 110592000000000 \approx 10^{14},
\ee
we expect the accuracy of the deductions in the table as long as the implementation does not have an error. In particular, we expect $d_{10} = 3$, which is compatible with the dimension conjecture. See also \cite{Hof15} and \cite{Sai17} for the preceding works for weights $\leq 9$.

\vspace{0.3in}
\addcontentsline{toc}{section}{Acknowledgements}
\noindent {\Large \bf Acknowledgements}
\vspace{0.2in}

\noindent
I thank S.\ Seki and M.\ Ono for introducing recent publication status of studies of FMZVs. I thank all people who helped me to learn mathematics and programming. I also thank my family.

%

\end{document}